\documentclass[a4paper, 12pt]{article}

\usepackage{amsmath, amssymb, amsthm}
\usepackage[dvips,final]{graphicx}
\usepackage{color}

\definecolor{blue3}{rgb}{0.05,0.05,0.5}
\definecolor{green1}{rgb}{0.2,0.6,0}
\definecolor{viol}{rgb}{.5,0,.3}
\definecolor{cyan2}{rgb}{0,.5,.5}

\numberwithin{equation}{section}

\newtheorem{theorem}{Theorem}[section]

\newtheorem{hypothesis}[theorem]{Hypothesis}

\usepackage{titlesec}
\titleformat{\section}{\bfseries}{\thesection}{1em}{}
\titleformat{\subsection}{\itshape}{\thesubsection}{1em}{}

\font\ctv msam10

\newcommand{\bbox}{\mbox{\ctv \symbol{4}}}
\def\QED{{$\hfill\bbox$}}
\newenvironment{pf}[1]{\par\vskip1mm{\noindent\it #1.}\ }{\QED\par\vskip2mm}

\def\bpf{\begin{pf}}
\def\epf{\end{pf}}

\newcommand{\real}{\mathbb{R}}
\newcommand{\R}{\mathbb{R}}
\newcommand{\nat}{\mathbb{N}}

\def\eps{\varepsilon}

\def\om{^{(m)}}

\def\for{\mbox{ for }\ }

\def\ale{\mbox{\ a.\,e. }}

\def\expe{\mathrm{e}}

\def\supess{\mathop{\mathrm{sup\,ess}}}

\def\mot#1{\left|#1\right|_{[0,t]}}

\def\Mot#1{\left\|#1\right\|_{[0,t]}}

\def\dehat{\hat}
\def\tht{\mathbf{u}}
\def\vtht{\mathbf{u}}

\def\play{\mathfrak{p}}

\def\dd{\mathrm{\,d}}

\def\PP{\mathcal{P}}
\def\UU{\mathcal{V}}

\def\htt{\hat\theta}
\def\PT{\PP(\theta)[q]}
\def\PTT{\PP_\theta(\theta)[q]}
\def\PTt{\PP_t(\theta)[q]}
\def\UT{\UU(\theta)[q]}
\def\UTT{\UU_\theta(\theta)[q]}
\def\UTt{\UU_t(\theta)[q]}
\def\PTh{\PP(\htt)[q]}
\def\PTTh{\PP_\theta(\htt)[q]}
\def\PTth{\PP_t(\htt)[q]}
\def\UTh{\UU(\htt)[q]}
\def\UTTh{\UU_\theta(\htt)[q]}
\def\UTth{\UU_t(\htt)[q]}
\def\PTm{\PP(\htt\om)[q\om]}
\def\PTTm{\PP_\theta(\htt\om)[q\om]}
\def\PTtm{\PP_t(\htt\om)[q\om]}
\def\UTm{\UU(\htt\om)[q\om]}
\def\UTTm{\UU_\theta(\htt\om)[q\om]}
\def\UTtm{\UU_t(\htt\om)[q\om]}
\def\PTTT{\PP_{\theta\theta}(\theta)[q]}
\def\UTTT{\UU_{\theta\theta}(\theta)[q]}
\def\PTTt{\PP_{\theta t}(\theta)[q]}
\def\UTTt{\UU_{\theta t}(\theta)[q]}

\def\ixr{\int_0^\infty\!\!\int_0^{\xi_r}}
\def\il{\int_0^\ell}

\def\be{\begin{equation}\label}
\def\ee{\end{equation}}

\begin{document}

\title{Oscillations of a temperature-dependent piezoelectric rod
\thanks{This research was supported by RVO: 67985840.}}

\author{Pavel~Krej\v{c}\'{i} and Giselle A. Monteiro
\thanks{Institute of Mathematics, Czech Academy of Sciences,
\v{Z}itn\'{a} 25, CZ-11567 Praha 1, Czech Republic, e-mail:~{\tt krejci@math.cas.cz},
{\tt gam@math.cas.cz}.}}

\maketitle

\begin{abstract}
Piezoelectricity of some materials has shown to have many applications, in particular in energy harvesting. 
Due to the inherent hysteresis in the characteristic of such materials, a number of hysteretic models have been proposed minding the fact that hysteresis losses may influence the efficiency of the process. However, hysteresis dissipation is accompanied with heat production, which in turn increases the temperature of the device and may change its physical characteristics. In this paper we propose a phenomenological model for electromechanical coupling in piezoelectric materials where temperature and feedback effects are taken into account. We prove the existence of solution for the resulting PDE system and show that the model is thermodynamically consistent. The main analytical tool is the inverse Preisach operator with temperature-dependent density.
\end{abstract}

\noindent
2010 AMS Classification: {34C55, 78A55, 49J15, 49J21}


\section*{Introduction}\label{int}

There are a lot of technological applications of multifunctional materials that spontaneously transform mechanical
energy into the electromagnetic one and vice versa. The best known examples are sensors and actuators for high
accuracy micropositioning or active damping of vibrations in real time, see \cite{ku,kku,rcl}.
Another important field of application is related to autonomous monitoring of bridges
and similar structures subject to permanent mechanical loading, see \cite{psg,xwsz}. 
A piezoelectric or magnetostrictive element is placed into the most exposed part of the construction and, by the effect
of mechanical vibrations, produces electric signal which is recorded and evaluated by an attached computer.
Simultaneously, the electromagnetic energy produced during the process is used for recharging the battery of the computer
and sending a wireless signal to the control point. Any anomalous behavior of the construction can therefore be
immediately detected and a possible problem can be fixed before irreversible damage occurs.

The main challenge is to model properly the electromechanical or magnetomechanical ``butterfly'' shaped curve.
A thermodynamic model for magnetostriction based on measurements of Galfenol carried out at the
University of Sannio at Benevento was proposed and analyzed
in \cite{dkv}. The underlying idea was motivated by the observation that both the magnetization hysteresis loops
and the magnetostrictive butterfly loops manifest a self-similar character parameterized by the applied stress.
This has led to the modeling hypothesis that all hysteresis phenomena can be described by one single Preisach operator
and its associated energy potential operator acting on an auxiliary self-similar variable. The butterfly-shaped
magnetostrictive curve then arises in a natural way from thermodynamic considerations involving the butterfly-shaped
Preisach potential operator. A similar model was shown to be applicable in piezoelectricity modeling in \cite{KK13}. Feedback effects have been taken into account in \cite{kk17}.

Here, we develop the idea of \cite{KK13} and include both the feedback and the temperature effects. More
specifically, we assume that the density function of the underlying Preisach operator depends on temperature.
The dependence cannot be arbitrary if we want to stay within the limits of the principles of thermodynamics.
We propose a formula for the free energy associated with the full thermo-electromechanical system
which is compatible with the Clausius-Duhem inequality. The 1D dynamics of an oscillating thermo-piezoelectric
rod is described by the mechanical momentum balance equation, the Gauss law, and the internal energy balance equation.
The resulting system is shown to admit a solution in an appropriate function space.

The main argument in the existence proof is a continuous inversion theorem for Preisach operators with temperature dependent density. This result has been presented in \cite{km1} and it extends considerably a similar statement in \cite{KK13}. 

The text is organized as follows. In Section \ref{mod} we present the model and show that the free energy can be chosen
so as to satisfy the principles of thermodynamics. An explicit formulation of the PDE system describing 1D oscillations
of a piezoelectric rod under thermal effects, the model hypotheses, and the statement of the existence result
are given in Section \ref{pde}. Section \ref{exi} is devoted to the proof of the main Theorem \ref{t1}. 
Results from \cite{km1} concerning the temperature-dependent Preisach inversion formula are collected in Section \ref{inv}.


\section{The model}\label{mod}

We consider the electric field $E$, the mechanical strain $\eps$, and the absolute temperature $\theta$
as state variables and the dielectric displacement
$D=D(\eps,E,\theta)$ as well as the mechanical stress $\sigma=\sigma(\eps,E,\theta)$
as state functions.

For the constitutive behavior of these quantities, we assume, similarly to \cite{dkv}, that
hysteresis effects are due to one single temperature-dependent Preisach operator $\PT$ with potential $\UT$
acting on an auxiliary state function $q = q(\eps,E)$, assuming that the Preisach density function
which determines the shape of the hysteresis loops depends on temperature.

Let us recall the definition of the Preisach model in the equivalent form of \cite{max}.
It is based on the concept of {\em play operator\/} $\play_r$ which is the mapping that with a given function
$q \in W^{1,1}(0,T)$ and a parameter $r>0$ associates the solution $\xi_r$ of the variational inequality
\be{vari}
\begin{array}{ll}
|q(t) - \xi_r(t)| \le r & \mbox{for all } t \in [0,T],\\
\dot \xi_r(t) (q(t) - \xi_r(t) - z) \ge 0 & \ale \ \forall z \in [-r,r],\\
\xi_r(0) = \max\{q(0) - r, \min\{0, q(0) + r\}\},
\end{array}
\ee
where the dot means the derivative with respect to $t$, and we denote $\xi_r(t) = \play_r[q](t)$.
As an immediate consequence of \eqref{vari}, we obtain the {\em hysteresis energy balance equation\/}
\be{eneq}
\dot\xi_r(t)(q(t) - \xi_r(t)) = r|\dot\xi_r(t)| \ \ale
\ee
This is indeed an energy balance, if we interpret $\dot\xi_r q$ as the power supplied to the system,
$\frac12 \xi_r^2$ as the potential energy, and $r|\dot\xi_r|$ as the dissipation rate.
The temperature-dependent Preisach operator $\PT$ is then defined by the integral
\be{prei}
\PT = \ixr \psi(\theta,r,v)\dd v\dd r,
\ee
where $\psi$ is a given nonnegative function called the {\em Preisach density\/}, which determines the shape
of the hysteresis loops and has to be determined experimentally, and
\be{pote}
\UT = \ixr v \psi(\theta,r,v)\dd v\dd r
\ee
is the associated Preisach potential. Note that for constant $\theta$, we have the Preisach energy inequality
\be{enerpr}
q \frac{\dd}{\dd t} \PT - \frac{\dd}{\dd t} \UT = \int_0^\infty r |\dot\xi_r| \psi(\theta,r,\xi_r)\dd r \ge 0.
\ee
As a temperature-dependent extension of the model in \cite{kk17}, we
assume that the polarization $P$ is given by the implicit relation
\be{pola}
P = \PT, \ q = \frac{1}{f(\eps)}(E - \alpha(\eps) P)
\ee
with a feedback parameter $\alpha(\eps)$, and a self-similarity function $f(\eps)$. We consider the stress $\sigma$,
the dielectric displacement $D$, and the free energy $F=F(\eps,E,\theta)$ of the form
\begin{eqnarray}
\sigma &=&\nu \eps_t + c \eps - e E
+ f'(\eps)\UT + \frac12 \alpha'(\eps) (\PT)^2 - \beta (\theta-\theta_c),
\label{sigma}\\ 
D&=& e \eps + \kappa E + \PT,
\label{diel}\\ 
F&=& F_0(\theta) + \frac{c}{2}\eps^2+\frac{\kappa}{2}E^2 + \frac{\gamma}{2} \eps_x^2 + f(\eps) \UT
+ \frac12 \alpha(\eps)(\PT)^2 - \beta(\theta - \theta_c)\eps,
\label{free}
\end{eqnarray}
where $\beta \ge 0$ is the thermal expansion coefficient, $\theta_c>0$ is a given reference temperature
(the room temperature, for example), $F_0(\theta)$ is the purely caloric part of the free energy which we specify later
and $q$ is given by \eqref{pola}. The term $\frac{\gamma}{2} \eps_x^2$ accounts for the couple stress,
see \cite{tou,falk}. The couple stress term is needed here in order to control the higher power terms in the
energy balance equation \eqref{PDE2} below, and we will comment on this issue later on.

We now check that the model is compatible with the principles of thermodynamics. The local energy balance equation
reads
\be{eba1}
U_t + Q_x = \sigma \eps_t + \gamma\eps_x\eps_{xt} + E D_t,
\ee
where $U = F+\theta S$ is the internal energy of the system and $Q$ is the heat flux that we assume
according to the Fourier law in the form
\be{four}
Q = -\mu \theta_x.
\ee
The interpretation of Eq.~\eqref{eba1} is the following. If we integrate \eqref{eba1} over a control interval $(a,b)$,
then the right-hand side represents the power supplied to the given volume, part of this power flows out of the interval
as the flux difference $Q(b) - Q(a)$, and the rest is used for the internal energy increase $\frac{\dd}{\dd t}\int_a^b U\dd x$.
The Second principle in Clausius-Duhem form states that there exists a state function called the entropy
such that 
\be{eba2}
S_t + \Big(\frac{Q}{\theta}\Big)_x \ge 0
\ee
for every process. Assuming the positivity of $\theta$ for the moment (and this will be proved later on), 
it is easy to see that if \eqref{eba1} holds, then \eqref{eba2} is satisfied provided the entropy
$S$ is chosen in such a way that
\be{cd}
\sigma \eps_t + E D_t + \gamma\eps_x\eps_{xt} - \theta_t S - F_t \ge 0
\ee
for every process. We claim that the right choice for the entropy is
\be{entr}
S = -F'_0(\theta) + \beta\eps + f(\eps) \big(q \PTT - \UTT\big),
\ee
where we denote
\be{preit}
\PTT = \ixr \psi_\theta(\theta,r,v)\dd v\dd t,\quad \UTT = \ixr v \psi_\theta(\theta,r,v)\dd v\dd t.
\ee
Indeed, a straightforward computation yields
\be{cd2}
\sigma \eps_t + E D_t + \gamma\eps_x\eps_{xt} - \theta_t S - F_t = \nu\eps_t^2
+ f(\eps)\left(\int_0^\infty r |(\xi_r)_t| \psi(\theta,r,\xi_r)\dd r\right) \ge 0
\ee
similarly as in \eqref{enerpr}.

We will consider a thermomechanical process in the domain $(x,t)\in (0,\ell)\times(0,T)$, where
$(0,\ell)$ is a space interval representing the 1D rod of length $\ell$ and $(0,T)$ is a given time interval.
Let $u(x,t)$ be the longitudinal displacement of a point $x$ of the rod at time $t$. Then
\be{strain}
\eps = u_x.
\ee
The full 1D system for unknown functions $u(x,t)$, $E(x,t)$, $\theta(x,t)$, $(x,t)\in (0,\ell)\times(0,T)$,
describing longitudinal oscillations of a thermo-piezoelectric rod consists of the energy balance \eqref{eba1}, 
and of the momentum balance and the Gauss law as follows:
\begin{equation}\label{pde0}
\begin{aligned}
\rho u_{tt} + \gamma u_{xxxx} - \sigma_x =&\ 0,\\
D_x=&\ 0.
\end{aligned}
\end{equation}
The equation $D_x = 0$ means that $D$ is a function of $t$ only, say, $D(x,t)=r(t)$,
where $r(t)$ is a function which is known from the boundary condition
$D(0,t)=D(\ell,t)=r(t)$, corresponding to an impressed (or measured) boundary current.
In order to simplify the analysis, we assume that $r(t)\equiv 0$, and prescribe also the simplest boundary conditions,
for the other unknowns, that is,
\begin{eqnarray}\label{bc1}
&& u=0 \ \mbox{ on }\ x=0,\ \gamma u_{xxx} - \sigma = 0\ \mbox{ on }\ x=\ell,\ u_{xx} = 0
\ \mbox{ on }\ x=0, \ell.\\ \label{bc2}
&&\theta_x = 0\ \mbox{ on }\ x=0, \ell.
\end{eqnarray}
The argument for more realistic non-homogeneous boundary conditions will be similar, just the formulas
would become a bit heavy. The condition $D=0$ means that
\begin{equation}\label{inverse}
e \eps +\kappa E + \PT = 0\,,
\end{equation}
and by virtue of \eqref{pola} we deduce the equation for $q$ in terms of $\eps$ and $\theta$
\be{q}
q + \frac{1+ \kappa\alpha(\eps)}{\kappa f(\eps)} \PT = - \frac{e\eps}{\kappa f(\eps)}.
\ee
It was shown in \cite{km1} that this equation determines $q$ uniquely in terms of $\theta$ and $\eps$,  
and the mapping $(\eps,\theta) \mapsto q$ is Lipschitz continuous with respect to the sup-norm (cf. Theorem \ref{p37}).


\section{Statement of the PDE problem}\label{pde}

Referring to \eqref{pola}, \eqref{sigma}, \eqref{free}, \eqref{entr}, \eqref{strain}, \eqref{inverse},
and putting $U = F + \theta S$, we rewrite the system \eqref{pde0} in variational form
\begin{eqnarray}\nonumber
\il (\rho u_{tt} w + \gamma u_{xx} w_{xx})\dd x &=& -\il \Big(\nu u_{xt} + c u_x + \frac{e}{\kappa} (e u_x +\PT)\\ \label{PDE1}
&&+ f'(u_x)\UT + \frac12 \alpha'(u_x) (\PT)^2 - \beta \theta\Big) w_x \dd x,\qquad\\ \nonumber
\il \big(-\theta F_0''(\theta)\theta_t\,z + \mu \theta_{x} z_x\big) \dd x &=&
\il \Big(\nu u_{xt}^2 - \beta\theta u_{xt} + f(u_x) (q \PTt - \UTt)\\ \label{PDE2}
&& - \theta \big(f(u_x)(q \PTT - \UTT)\big)_t\Big) z \dd x,
\end{eqnarray}
for every test functions $w \in W^{2,2}(0,\ell) \cap W^{1,2}_0(0,\ell)$ and $z \in W^{1,2}(0,\ell)$, where we denote
\be{w120}
W^{1,2}_0(0,\ell) = \{v \in W^{1,2}(0,\ell): v(0) = 0\},
\ee
and
\be{ptt}
\PTt = \int_0^\infty (\xi_r)_t \psi(\theta,r,\xi_r)\dd r, \ \UTt = \int_0^\infty (\xi_r)_t \xi_r \psi(\theta,r,\xi_r)\dd r
\ee
with $q$ defined as the solution of \eqref{q} and $\xi_r$ as in \eqref{vari}. The term
\be{dis}
f(u_x) (q \PTt - \UTt) = f(u_x) \int_0^\infty r |(\xi_r)_t| \psi(\theta,r,\xi_r)\dd r \ge 0
\ee
is the hysteresis dissipation rate as part of the entropy production rate in \eqref{cd2} and appears in the energy balance
\eqref{PDE2} as heat source.

The function $c_V(\theta) = -\theta F_0''(\theta)$ is the specific heat capacity. For example, the choice
$F_0(\theta) = - c_0\theta \log(\theta/\theta_c)$ would correspond to the assumption that $c_V(\theta) = c_0$
is constant.

We now check that the total energy of the system is formally conserved during the evolution. 
We test Eq.~\eqref{PDE1} by $w=u_t$, Eq.~\eqref{PDE2} by $z=1$, and sum up. We obtain
\begin{eqnarray} \nonumber
&&\frac{\dd}{\dd t} \int_0^\ell \Big(F_0(\theta) -\theta F_0'(\theta) + \frac{\rho}{2} u_t^2 + \frac{c}{2} u_x^2
+ \frac{\gamma}{2} u_{xx}^2 + f(u_x)\UT
+ \frac12 \alpha(u_x)(\PT)^2\\ \label{cons}
&&\quad + \frac{1}{2\kappa}\big(e u_x +\PT\big)^2 + \theta f(u_x)(q \PTT - \UTT)\Big) \dd x = 0.
\end{eqnarray}
Indeed, the expression in \eqref{cons} under the time derivative is the total energy of the system, and its
time derivative is zero.

We now formulate the hypotheses that are assumed to hold. For practical reasons, we list separately
the assumptions about the non-hysteretic terms in Problem \eqref{PDE1}--\eqref{PDE2} (Hypothesis \ref{h1})
and about the Preisach operator \eqref{prei} (Hypothesis \ref{h2}).

\begin{hypothesis}\label{h1}
The functions occurring in the statement of the problem fulfill the conditions
\begin{itemize}
\item[{\rm (i)}] The function $f:\real \to \real$ is bounded from above and from below by constants
$0< f_0 \le f(\eps) \le f_1$, both $f$ and $f'$ are Lipschitz continuous, and the function
$\eps \mapsto (1+|\eps| |f'(\eps)|)$ is bounded;
\item[{\rm (ii)}] The function $\alpha:\real \to \real$ is bounded, 
$1+\kappa\alpha(\eps) \ge 0$, and both $\alpha$ and $\alpha'$ are Lipschitz continuous, and we put
$A^*:= \sup_{\eps \in \real}\frac{1+\kappa \alpha(\eps)}{\kappa f(\eps)}$;
\item[{\rm (iii)}] The function $\theta\mapsto C_V(\theta) := -\theta F_0''(\theta)$ is continuous,
and there exists a constant $c_0 > 0$ such that $C_V(\theta) \ge c_0 (1 + (\theta^+)^{1/3})$ for all $\theta \in \real$;
\item[{\rm (iv)}] The initial conditions have the regularity $u^0 \in W^{4,2}(0,\ell)\cap W^{1,2}_0(0,\ell)$,
$u^0_{xx}(0) = u^0_{xx}(\ell) = 0$, $u^1 \in W^{2,2}(0,\ell)\cap W^{1,2}_0(0,\ell)$, $\theta^0 \in W^{1,2}(0,\ell)$,
$\theta^0(x) > 0$ \ale
\end{itemize}
\end{hypothesis}

\begin{hypothesis}\label{h2}
The Preisach density $\psi(\theta,r,v) \in L^1(\real\times (0,\infty) \times \real)$ in \eqref{prei} is such that
the function $\theta \mapsto \psi(\theta,r,v)$ is of class $C^2$ in $\real$ for a.\,e. $(r,v)\in (0,\infty)\times \real$.
Let
$$
g(\theta,r,v) = \int_0^v \psi(\theta,r,v')\dd v'.
$$
We assume that there exists constants $\Psi_0> 0$ and $\delta \in (0,1)$ such that
\begin{itemize}
\item[{\rm (i)}] $\psi(\theta,r,v) = \psi(\theta,r,-v) \ge 0 \ale$;
\item[{\rm (ii)}] $\psi_\theta(\theta,r,v) = 0\ \ale\ \for \theta\le 0$;
\item[{\rm (iii)}] $\int_0^\infty(1+r)\psi(\theta,r,v)\dd r \le \Psi_0 \ \ \forall \theta\in \real, v\in \real$;
\item[{\rm (iv)}] $\int_0^\infty\int_0^\infty \psi(\theta,r,v)\dd v \dd r \le \Psi_0
\ \ \forall \theta\in \real$;
\item[{\rm (v)}] $\int_0^\infty\int_0^\infty v \psi(\theta,r,v)\dd v \dd r \le \Psi_0 (1+\theta^+) \ \ \forall \theta\in \real$;
\item[{\rm (vi)}] $\int_0^\infty |v|\, \psi(\theta,r,v) \dd r \le \Psi_0 (1+(\theta^+)^{1/6})
\ \ \forall \theta\in \real, v\in \real$;
\item[{\rm (vii)}] $\int_0^\infty\int_0^\infty (v+\theta(1+r)) |\psi_\theta(\theta,r,v)| \dd v\dd r \le \Psi_0\ \ \forall \theta\in \real$;
\item[{\rm (viii)}] $\int_0^\infty (1+\theta)(1+r) |\psi_\theta(\theta,r,v)| \dd r \le \Psi_0
\ \ \forall \theta>0, v \in \real$;
\item[{\rm (ix)}] $\left|\int_0^\infty\int_0^K (1+\theta) g_{\theta}(\theta,r,v)\dd v \dd r\right| \le \Psi_0
\ \ \forall \theta>0, K>0$;
\item[{\rm (x)}] $\int_0^\infty\int_0^\infty |\psi_{\theta}(\theta,r,v)|\dd v \dd r
\le \delta/f_1 \ \ \forall \theta>0$;
\item[{\rm (xi)}] $\left|\int_0^\infty\int_0^\infty r(1+\theta) \psi_{\theta\theta}(\theta,r,v)\dd v \dd r\right|
\le \delta/f_1 \ \ \forall \theta>0$;
\item[{\rm (xii)}] $\left|\int_0^\infty\int_0^K (1+\theta) g_{\theta\theta}(\theta,r,v)\dd v \dd r\right| \le \delta/f_1
\ \ \forall \theta>0, K>0$;
\item[{\rm (xiii)}] $(7 + 3A^*\Psi_0)\delta \le c_0/2$ with $A^*$ as in Hypothesis \ref{h1}\,(ii).
\end{itemize}
\end{hypothesis}

A typical function $\psi$ satisfying Hypothesis \ref{h2} might have the form
\be{psi}
\psi(\theta,r,v) =
\left\{
\begin{array}{ll}
a(r) \phi(|v| - h(\log(1{+}\theta))) & \for \ \theta >0,\\
a(r) \phi(|v|) & \for \ \theta\le 0,
\end{array}
\right.
\ee
where $a\ge 0$ is a function in $L^1(0,\infty)$ such that $r \mapsto r a(r)$ belongs to $L^1(0,\infty)$,
$\phi\ge 0$ is a $C^2$-function with compact support in the interval $(0,1)$,
and $h:[0,\infty) \to [0,\infty)$ is a $C^2$-function such that $h(0) = h'(0) = h''(0) = 0$, 
$0 \le h'(s) \le 1$, $|h''(s)| \le 1$ for all $s\ge 0$, for example $h(s) = \frac{h_0 s^3}{1+s^2}$.

To see that Hypothesis \ref{h2} is fulfilled with the choice \eqref{psi} of $\psi$, note that for $\theta>0$ and $v\ge 0$ we have
\begin{eqnarray*}
\psi_\theta(\theta,r,v) &=& -a(r) \frac{h'(\log(1{+}\theta))}{1{+}\theta} \phi'(v - h(\log(1{+}\theta))),\\
\psi_{\theta\theta}(\theta,r,v) &=& a(r) \left(\frac{h'(\log(1{+}\theta))}{1{+}\theta}\right)^2 \phi''(v - h(\log(1{+}\theta)))\\
&&+\, a(r) \frac{h'(\log(1{+}\theta)){-} h''(\log(1{+}\theta))}{(1{+}\theta)^2}\phi'(v - h(\log(1{+}\theta))),\\
g_\theta(\theta,r,v) &=& -a(r) \frac{h'(\log(1{+}\theta))}{1{+}\theta} \phi(v - h(\log(1{+}\theta))),\\
g_{\theta\theta}(\theta,r,v) &=& a(r)\left(\frac{h'(\log(1{+}\theta))}{1{+}\theta}\right)^2 \phi'(v - h(\log(1{+}\theta)))\\
&&+\, a(r) \frac{h'(\log(1{+}\theta)){-} h''(\log(1{+}\theta))}{(1{+}\theta)^2}\phi(v - h(\log(1{+}\theta))),\\
\int_0^\infty v\psi(\theta,r,v)\dd v &=& a(r) \int_0^1 (s + h(\log(1{+}\theta)))\phi(s)\dd s,
\end{eqnarray*}
and it suffices to impose suitable assumptions on the $L^1$-norm of the function $r\mapsto (1+r)a(r)$.

We now show formally how Hypothesis \ref{h2} will be used to estimate the hysteresis terms on the right-hand sides of
\eqref{PDE1}--\eqref{PDE2}. We will use the following splitting:
$$
\begin{aligned}
q\PTTT - \UTTT =& \ixr(q-v) \psi_{\theta\theta}(\theta,r,v)\dd v\dd r\\
=& \ixr ((q-\xi_r) + (\xi_r-v)) \psi_{\theta\theta}(\theta,r,v)\dd v\dd r,
\end{aligned}
$$
where
$$
\left|\ixr(q-\xi_r) \psi_{\theta\theta}(\theta,r,v)\dd v\dd r\right|
\le \ixr r \left|\psi_{\theta\theta}(\theta,r,v)\right|\dd v\dd r\,,
$$
and
$$
\left|\ixr (\xi_r-v) \psi_{\theta\theta}(\theta,r,v)\dd v\dd r\right|
= \left|\ixr g_{\theta\theta}(\theta,r,v)\dd v\dd r\right|\,,
$$
hence
\be{he}
|q\PTTT - \UTTT| \le \ixr \left|\psi_{\theta\theta}(\theta,r,v)\right|\dd v\dd r
+ \left|\ixr g_{\theta\theta}(\theta,r,v)\dd v\dd r\right|,
\ee
and similarly
\be{he0}
|q\PTT - \UTT| \le \ixr r \left|\psi_{\theta}(\theta,r,v)\right|\dd v\dd r
+ \left|\ixr g_{\theta}(\theta,r,v)\dd v\dd r\right|.
\ee
In the series of inequalities below, the superscript (iii)--(xii) refers to the
corresponding item in Hypothesis \ref{h2}, and $C$ denotes any positive constant depending only on the data
of the problem.
\begin{eqnarray} \label{he1}
|\PT| &=& \left|\ixr \psi(\theta,r,v)\dd v\dd r\right| \stackrel{(iv)}{\le} \Psi_0,\\ \label{he2}
0 \le \UT &=& \left|\ixr v \psi(\theta,r,v)\dd v\dd r\right| \stackrel{(v)}{\le} \Psi_0 (1+\theta^+),\\ \label{he3}
|\PTt| &=& \left|\int_0^\infty (\xi_r)_t \psi(\theta,r,\xi_r)\dd r\right| \stackrel{(iii)}{\le} \Psi_0 |q_t|,\\ 
\label{he4}
|\theta_t\PTT| &=& |\theta_t|\left|\ixr \psi_\theta(\theta,r,v)\dd v\dd r\right|\stackrel{(x)}{\le} |\theta_t|\delta/f_1,\\ 
\nonumber
|\UTt| &=& \left|\int_0^\infty (\xi_r)_t \xi_r \psi(\theta,r,\xi_r)\dd r\right|
\le |q_t| \int_0^\infty |\xi_r| \psi(\theta,r,\xi_r)\dd r\\ \label{he5}
&\stackrel{(vi)}{\le}& \Psi_0 |q_t| (1+(\theta^+)^{1/6}),\\ 
\label{he6}
0 \le f(u_x)(q\PTt - \UTt) &=& f(u_x)\int_0^\infty r(\xi_r)_t \psi(\theta,r,\xi_r)\dd r
\stackrel{(iii)}{\le} \Psi_0 f(u_x)|q_t|,\\ 
\nonumber
|\theta (q\PTT {-} \UTT)| &\stackrel{\eqref{he0}}{\le}
& C \theta\Big(\ixr r \big|\psi_{\theta}(\theta,r,v)\big|\dd v\dd r\\ 
\label{he7}
&&+\, \Big|\ixr g_{\theta}(\theta,r,v)\dd v\dd r\Big|\Big)\stackrel{(vii),(ix)}{\le} C,\\ 
\nonumber
\theta f(u_x)|q_t\PTT| &=& \theta f(u_x) |q_t| \left|\ixr \psi_\theta(\theta,r,v)\dd v\dd r\right|\\ \label{he8}
&\stackrel{(vii)}{\le}& \Psi_0 f(u_x)|q_t|,\\ \label{he9}
f(u_x)|\theta(q\PTTt - \UTTt)| &=& f(u_x)\theta\left|\int_0^\infty(\xi_r)_t \psi_\theta(\theta,r,\xi_r)\dd r\right|
\stackrel{(viii)}{\le} \Psi_0 f(u_x)|q_t|,\\ 
\nonumber
|f(u_x)\theta_t \theta (q\PTTT {-} \UTTT)| &\stackrel{\eqref{he}}{\le}&
f_1 |\theta_t| \theta\Big(\ixr r \big|\psi_{\theta\theta}(\theta,r,v)\big|\dd v\dd r\\ \label{he10}
&&+\, \Big|\ixr g_{\theta\theta}(\theta,r,v)\dd v\dd r\Big|\Big)\stackrel{(xi),(xii)}{\le} 2\delta|\theta_t|.
\end{eqnarray}
Eq.~\eqref{q} can be written in the form
\be{q1}
q + A(u_x) \PT = B(u_x)
\ee
with Lipschitz continuous functions $A, B$. Moreover, by Hypothesis \ref{h1}\,(ii), we have
\be{q2}
A(u_x) \le A^*:= \sup_{\eps \in \real}\frac{1+\kappa \alpha(\eps)}{\kappa f(\eps)}. 
\ee
Differentiating \eqref{q1} in $t$ gives
\be{q3}
q_t + A'(u_x) u_{xt} \PT + A(u_x) (\PTt + \theta_t \PTT) = B'(u_x) u_{xt}.
\ee
Let $C$ denote here again and in the sequel any constant depending only on the data of the problem.
Using the fact that $q_t\PTt \ge 0$ a.~e., we obtain
\be{q4}
|q_t| \le C |u_{xt}| + \frac{A^*\delta}{f_1} |\theta_t| \ \ale
\ee
We are now ready to state the main existence result.

\begin{theorem}\label{t1}
Let Hypotheses \ref{h1}, \ref{h2} hold, and let $T>0$ be a given final time. Then system \eqref{PDE1}--\eqref{PDE2}
with initial conditions
\be{ini}
u(x,0) = u^0(x), \ u_t(x,0) = u^1(x), \ \theta(x,0) = \theta^0(x) \ \for x \in (0,\ell)
\ee
admits a solution with the regularity
$$
u_{xxt}, \theta_x \in L^\infty(0,T;L^2(0,\ell)), \ u_{xtt}, \theta_t, \theta_{xx} \in L^2((0,\ell)\times(0,T)),
$$
$\theta(x,t) \ge 0$ for all $(x,t) \in (0,\ell)\times(0,T)$.
\end{theorem}


\section{Proof of the existence theorem}\label{exi}

Theorem \ref{t1} will be proved in several steps. We first fix a cut-off parameter $R>0$, define a cut-off mapping
$K_R(z) = \min\{R, z^+\}$ for $z \in \real$, replace $\theta$ at critical places with $\htt = K_R(\theta)$
and consider instead of \eqref{PDE1}--\eqref{PDE2} the truncated system
\begin{eqnarray}\nonumber
\il (\rho u_{tt} w + \gamma u_{xx} w_{xx})\dd x &=& -\il \Big(\nu u_{xt} + c u_x + \frac{e}{\kappa} (e u_x +\PTh)\\ \label{PDE3}
&&+ f'(u_x)\UTh + \frac12 \alpha'(u_x) (\PTh)^2 - \beta \htt \Big) w_x \dd x,\qquad\\ \nonumber
\il \left(C_V(\htt)\theta_t\,z + \mu \theta_{x}z_x\right)\dd x &=&
\il\Big(\nu K_R(u_{xt}^2) - \beta \htt u_{xt} + f(u_x) (q \PTth - \UTth)\\ \label{PDE4}
&&- \htt \big(f(u_x)(q \PTTh - \UTTh)\big)_t\Big)z \dd x
\end{eqnarray}
for all test functions $w \in W^{2,2}(0,\ell) \cap W^{1,2}_0(0,\ell)$ and $z \in W^{1,2}(0,\ell)$,
with $q$ defined as the solution of the equation
\be{qh}
q + \frac{1+ \kappa\alpha(u_x)}{\kappa f(u_x)} \PTh = - \frac{e u_x}{\kappa f(u_x)}.
\ee


\subsection{ODE approximation}\label{odes}

System \eqref{PDE3}--\eqref{PDE4} will be solved by Galerkin approximations. We choose orthonormal bases in $L^2(0,\ell)$
\be{ortho}
s_k(x) = \sqrt{\frac{2}{\ell}}\sin \frac{k\pi}{\ell} x, \ c_k(x) = \sqrt{\frac{2}{\ell}}\cos \frac{k\pi}{\ell} x
\ \for k\in \nat, \ c_0(x) = \frac{1}{\sqrt{\ell}}.
\ee
For each $m \in \nat$ we determine the coefficients $u_k(t), \theta_k(t)$ in the expansions
\be{gale}
u\om(x,t) = \sum_{k=1}^m u_k(t) s_k(x), \ \theta\om(x,t) = \sum_{k=0}^m \theta_k(t) c_k(x)
\ee
as solutions of the ODE system
\begin{eqnarray}\nonumber
&&\il (\rho u\om_{tt} s_k + \gamma u\om_{xx} s_k'')\dd x = -\il \Big(\nu u\om_{xt} + c u\om_x
+ \frac{e}{\kappa} (e u\om_x +\PTm)\\ \label{PDE3m}
&&\quad + f'(u\om_x)\UTm + \frac12 \alpha'(u\om_x) (\PTm)^2 - \beta \htt\om\Big) s_k' \dd x,\\ \nonumber
&&\il \left(\Big(\frac{1}{m} \theta\om_{tt} + C_V(\htt\om)\theta\om_t\Big) c_k + \mu \theta\om_{x} c_k'\right)\dd x =
\il \Big(\nu K_R((u\om_{xt})^2)\\ \nonumber
&&\quad - \beta \htt\om u\om_{xt} + f(u\om_x) (q\om \PTtm - \UTtm)\\[2mm] \label{PDE4m}
&&\quad - \htt\om \big(f(u\om_x)(q\om \PTTm - \UTTm)\big)_t\Big) c_k \dd x
\end{eqnarray}
with initial conditions
\be{inim}
\begin{split}
& u_k(0) = \il u^0(x) s_k(x) \dd x, \qquad \dot u_k(0) = \il u^1(x) s_k(x) \dd x,
\\
&\theta_k(0) = \il \theta^0(x) c_k(x) \dd x, \dot\theta_k(0) = 0,
\end{split}
\ee
for $k=0, \dots, m$ with the convention $s_0(x) = 0$. In \eqref{PDE3m}--\eqref{inim},
the prime denotes differentiation with respect to $x$, the dot denotes as before differentiation with respect to $t$, 
we use the notation
\be{htt}
\htt\om = K_R(\theta\om),
\ee
and $q\om$ is defined as the solution of the equation
\begin{equation*}
q\om + \frac{1+ \kappa\alpha(u\om_x)}{\kappa f(u\om_x)} \PP(\theta\om)[q\om] = - \frac{e u\om_x}{\kappa f(u\om_x)}.
\end{equation*}
For each $m\in\nat$, this is an ODE system with locally Lipschitz right hand side, hence the initial value problem
\eqref{PDE3m}-\eqref{inim} admits
a unique local solution in a maximal interval $[0, T_m)$. We now derive a series of estimates in $[0,T_m)$
which will imply that $T_m = T$ and that the sequence $(u\om, \theta\om)$ contains a convergent subsequence
which converges to a solution of \eqref{PDE3}--\eqref{PDE4}.
 We 
denote by $C_R$ any constant depending only on the data and
on the cut-off parameter $R$, but independent of $m$, and by $C$ as before any constant depending possibly on the data
and independent of $R$ and $m$.

We start by multiplying \eqref{PDE3m} by $\dot u_k(t)$, summing over $k = 1, \dots, m$, and integrating with respect to $t$ from $0$ to some $\tau\in[0, T_m)$. The temperature-dependent terms on the right-hand side of \eqref{PDE3m} are bounded by $C (1+\htt\om)$.
Moreover, the regularity of the initial conditions \eqref{inim} ensures that 
\begin{align*}
\il\left(\rho|u\om_t|^2 {+} \big(c+\frac{e}{\kappa}\big)|u\om_x|^2 {+}\gamma |u\om_{xx}|^2\right)(x,0) \dd x 
\le C.
\end{align*}
As a consequence we obtain
\be{es0}
\il\left(|u\om_t|^2 {+} |u\om_x|^2 {+} |u\om_{xx}|^2\right)(x,\tau) \dd x {+} \int_0^\tau\!\!\il |u\om_{xt}|^2 \dd x\dd t
\le C\left(1 {+} \int_0^\tau\!\!\il |\htt\om|^2 \dd x\dd t\right).
\ee
In particular, since $|\htt\om| \le R$, we have
\be{es1}
\il\left(|u\om_t|^2 + |u\om_x|^2 + |u\om_{xx}|^2\right)(x,\tau) \dd x + \int_0^\tau\il |u\om_{xt}|^2 \dd x\dd t \le C_R. 
\ee
We further multiply \eqref{PDE4m} by $\dot\theta_k(t)$,
sum over $k = 0, \dots, m$ and integrate with respect to $t$ from $0$ to $\tau$. This yields
\be{es2}
\begin{aligned}
&\hspace{-10mm}\il\left(\frac1{m}|\theta\om_t|^2 + \mu |\theta\om_x|^2\right)(x,\tau) \dd x
+ c_0\int_0^\tau\il (1+(\htt\om)^{1/3})|\theta\om_{t}|^2 \dd x\dd t\\
\le&\ \il \mu |\theta\om_x|^2(x,0) \dd x + \int_0^\tau\il \sum_{i=1}^4 H_i(x,t) \theta\om_t\dd x\dd t,
\end{aligned}
\ee
where we have
$$
\begin{aligned}
H_1(x,t) =&\ \nu K_R((u\om_{xt})^2),\\
H_2(x,t) =&\ - \beta \htt\om u\om_{xt},\\
H_3(x,t) =&\ f(u\om_x) (q\om \PTtm - \UTtm),\\
H_4(x,t) =&\ - \htt\om \big(f(u\om_x)(q\om \PTTm - \UTTm)\big)_t.
\end{aligned}
$$
Using \eqref{he6} and \eqref{q4} we get
$$
\begin{aligned}
\left|\int_0^\tau\il H_1(x,t) \theta\om_t\dd x\dd t\right|
\le &\ C \int_0^\tau\il K_R^2((u\om_{xt})^2)\dd x\dd t + \delta \int_0^\tau\il |\theta\om_{t}|^2 \dd x\dd t,\\
\left|\int_0^\tau\il H_2(x,t) \theta\om_t\dd x\dd t\right|
\le &\ C \int_0^\tau\il (\htt\om)^2 |u\om_{xt}|^2\dd x\dd t + \delta \int_0^\tau\il |\theta\om_{t}|^2 \dd x\dd t,\\
\left|\int_0^\tau\il H_3(x,t) \theta\om_t\dd x\dd t\right|
\le &\ C \int_0^\tau\il |u\om_{xt}| |\theta\om_t|\dd x\dd t
+ \delta A^*\Psi_0 \int_0^\tau\il |\theta\om_{t}|^2 \dd x\dd t\\
\le &\ C \int_0^\tau\il |u\om_{xt}|^2 \dd x\dd t
+ \delta (1+ A^*\Psi_0) \int_0^\tau\il |\theta\om_{t}|^2 \dd x\dd t,
\end{aligned}
$$
while the estimates \eqref{he7}--\eqref{he10}, \eqref{q4} imply
$$
\begin{aligned}
&\left|\int_0^\tau\il H_4(x,t) \theta\om_t\dd x\dd t\right|
\\&\ 
\le C \int_0^\tau\il |u\om_{xt}| |\theta\om_t|\dd x\dd t
+ 2\,\Psi_0\,f_1 \int_0^\tau\il |q\om_t||\theta\om_{t}|+2\,\delta\int_0^\tau\il |\htt\om_{t}||\theta\om_{t}| \dd x\dd t
\\&\ 
\le  C \int_0^\tau\il |u\om_{xt}| |\theta\om_t|\dd x\dd t
+ 2\delta (1+ A^*\Psi_0) \int_0^\tau\il |\theta\om_{t}|^2 \dd x\dd t
\\&\ 
\le  C \int_0^\tau\il |u\om_{xt}|^2 \dd x\dd t
+ 2\delta (2+ A^*\Psi_0) \int_0^\tau\il |\theta\om_{t}|^2 \dd x\dd t.
\end{aligned}
$$
Moreover, the regularity of the initial conditions \eqref{inim} guarantees that $\il |\theta\om_x|^2(x,0) \dd x\le C$. Recall that by Hypothesis \ref{h2}\,(xiii) we have $(7 + 3A^*\Psi_0)\delta \le c_0/2$, thus the computations above together with \eqref{es2} yield
\be{es2a}
\begin{aligned}
&\hspace{-10mm}\il\left(\frac1{m}|\theta\om_t|^2 + |\theta\om_x|^2\right)(x,\tau) \dd x
+ \int_0^\tau\il (1+(\htt\om)^{1/3})|\theta\om_{t}|^2 \dd x\dd t\\
\le&\ C\left(1+ \int_0^\tau\il \left(K_R^2((u\om_{xt})^2) + (1+(\htt\om)^2)|u\om_{xt}|^2\right)\dd x\dd t\right) \le C_R.
\end{aligned}
\ee
The next estimate is obtained by differentiating \eqref{PDE3m} in $t$, multiplying by $\ddot u_k(t)$, 
summing over $k = 1, \dots, m$, and integrating with respect to $t$ from $0$ to some $\tau$. This leads to 
\begin{equation*}
\begin{aligned}
&\hspace{-10mm}\il\left(|u\om_{tt}|^2 + |u\om_{xt}|^2 + |u\om_{xxt}|^2\right)(x,\tau) \dd x
+ \int_0^\tau\!\!\il |u\om_{xtt}|^2 \dd x\dd t\\
\le&\ C\left(1{+}\int_0^\tau\!\!\il |\htt\om_t|^2 \dd x\dd t + \sum_{i=1}^5\int_0^\tau\il \big(G_i(x,t)\big)^2\dd x\dd t\right),
\end{aligned}
\end{equation*}
where $G_i(x,t)$ are given by
$$
\begin{aligned}
G_1(x,t) =&\ \frac{e}{\kappa}\,(\PTm)_t,\\
G_2(x,t) =&\ - f''(u\om_x)u\om_{xt}\UTm,\\
G_3(x,t) =&\ f'(u\om_x) (\UTm)_t,\\
G_4(x,t) =&\ \frac{1}{2}\,\alpha''(u\om_x)u\om_{xt}(\PTm)^2\\
G_5(x,t) =&\ \alpha'(u\om_x)\PTm\,(\PTm)_t.
\end{aligned}
$$
Using \eqref{he1}--\eqref{he5}, \eqref{q4}, we can show that $C (|u\om_{xt}| + |\theta\om_t|)$ is an upper bound for $G_1,G_4$ and $G_5$, as well as 
$$
|G_2(x,t)| \le\ C\, |u\om_{xt}|\,(1 + \htt\om\,),\qquad
|G_3(x,t)| \le\ C (|u\om_{xt}| + |\theta\om_t|) (1+(\htt\om)^{1/6}\,).
$$ 
Therefore,
\begin{eqnarray}\nonumber
&&\hspace{-10mm}\il\left(|u\om_{tt}|^2 + |u\om_{xt}|^2 + |u\om_{xxt}|^2\right)(x,\tau) \dd x
+ \int_0^\tau\!\!\il |u\om_{xtt}|^2 \dd x\dd t \\ \label{es3}
&\le& C\left(1 {+}\int_0^\tau\!\!\il \big((1+(\htt\om)^2)|u\om_{xt}|^2
+ (1+(\htt\om)^{1/3})|\theta\om_t|^2\big) \dd x\dd t\right).
\end{eqnarray}
Combining \eqref{es2a} with \eqref{es3} we thus conclude that
\begin{eqnarray}\nonumber
&&\hspace{-10mm}\il\left(\frac1{m}|\theta\om_t|^2 + |\theta\om_x|^2 + |u\om_{tt}|^2 + |u\om_{xt}|^2
+ |u\om_{xxt}|^2\right)(x,\tau) \dd x\\ \nonumber
&&+ \int_0^\tau\!\!\il \left((1+(\htt\om)^{1/3})|\theta\om_{t}|^2 + |u\om_{xtt}|^2\right) \dd x\dd t \\ \label{es3a}
&\le& C\left(1 {+}\int_0^\tau\!\!\il \big(K_R^2((u\om_{xt})^2) + (1+(\htt\om)^2)|u\om_{xt}|^2\big) \dd x\dd t\right)
\le C_R,
\end{eqnarray}
for every $m\in\nat$ and $\tau\in[0,T_m)$.
Keeping $R$ fixed, by letting $m\to \infty$, and using compact embedding formulas, we can find
functions $(u,\theta)$ and a subsequence of $(u\om, \theta\om)$ (still indexed by $m$ for simplicity) such that
$$
\begin{array}{ll}
u\om \to u & \mbox{strongly in }\ C([0,\ell]\times [0,T])\\[2mm]
u\om_x \to u_x & \mbox{strongly in }\ C([0,\ell]\times [0,T])\\[2mm]
u\om_t \to u_t & \mbox{strongly in }\ C([0,\ell]\times [0,T])\\[2mm]
u\om_{xt} \to u_{xt} & \mbox{strongly in }\ C([0,\ell]\times [0,T])\\[2mm]
u\om_{tt} \to u_{tt} & \mbox{weakly-star in }\ L^\infty(0,T;L^2(0,\ell))\\[2mm]
u\om_{xx} \to u_{xx} & \mbox{weakly-star in }\ L^\infty(0,T;L^2(0,\ell))\\[2mm]
u\om_{xxt} \to u_{xxt} & \mbox{weakly-star in }\ L^\infty(0,T;L^2(0,\ell))\\[2mm]
u\om_{xtt} \to u_{xtt} & \mbox{weakly in }\ L^2((0,\ell)\times (0,T))\\[2mm]
\theta\om \to \theta & \mbox{strongly in }\ C([0,\ell]\times [0,T])\\[2mm]
\theta\om_{t} \to \theta_{t} & \mbox{weakly in }\ L^2((0,\ell)\times (0,T))\\[2mm]
\theta\om_{x} \to \theta_{x} & \mbox{weakly-star in }\ L^\infty(0,T;L^2(0,\ell)).
\end{array}
$$
We can therefore pass to the limit in \eqref{PDE3m}--\eqref{PDE4m} as $m \to \infty$ and conclude
that $(u,\theta)$ is a solution of \eqref{PDE3}--\eqref{PDE4} with the regularity as in Theorem \ref{t1}.
Indeed, the $L^2$-regularity of $\theta_{xx}$ is obtained by comparison with the other terms in \eqref{PDE4}. 
Moreover, since the terms on the right hand sides of \eqref{es3a} and \eqref{es0} converge strongly, we can pass to the limit
and conclude that the solution to \eqref{PDE3}--\eqref{PDE4} satisfies for every $\tau\in [0,T]$ the estimates
\begin{eqnarray}\nonumber
&&\hspace{-10mm}\il\left(|\theta_x|^2 + |u_{tt}|^2 + |u_{xt}|^2 + |u_{xxt}|^2\right)(x,\tau) \dd x
+ \int_0^\tau\!\!\il \left((1+\htt^{1/3})|\theta_{t}|^2 + |u_{xtt}|^2\right) \dd x\dd t \\ \label{es3b}
&\le& C\left(1 {+}\int_0^\tau\!\!\il \big(K_R^2((u_{xt})^2) + (1+\htt^2)|u_{xt}|^2\big) \dd x\dd t\right),\\ \label{es3c}
&&\hspace{-10mm}\il\left(|u_t|^2 {+} |u_x|^2 {+} |u_{xx}|^2\right)(x,\tau) \dd x
+ \int_0^\tau\!\!\il |u_{xt}|^2 \dd x\dd t \le C\left(1 {+} \int_0^\tau\!\!\il |\htt|^2 \dd x\dd t\right). 
\end{eqnarray}
with a constant $C$ independent of $R$.


\subsection{Positivity of temperature}\label{post}

Firstly, note that $\htt\theta^-=0$. Hence, by testing \eqref{PDE4} by $z=-\theta^-$ we get
\begin{equation*}
-\il \Big(C_V(\htt)\theta_t \theta^- + \mu \theta_x (\theta^-)_x\Big)\dd x 
	= -\il \Big(\nu K_R(u_{xt}^2)\theta^- + f(u_x) (q \PTth - \UTth)\theta^-\Big)\dd x.
\end{equation*}
Since the two terms on the right-hand side are non-negative, cf.~\eqref{dis}, using the fact that
$C_V(K_R(\theta))\theta^- = C_V(0)\theta^-$, we obtain
\be{pos1}
-C_V(0)\il \theta_t\theta^- \dd x - \il \mu \theta_x (\theta^-)_x \dd x \le 0, 
\ee
that is,
\be{pos2}
\frac{C_V(0)}{2}\frac{\dd}{\dd t} \il |\theta^-|^2\dd x + \il \mu \big|(\theta^-)_x\big|^2 \dd x \le 0.
\ee
By Hypothesis \ref{h1}\,(v) we have $\theta^-(x,0) = 0$, hence $\theta^-(x,t) = 0$ for all $(x,t) \in (0,\ell)\times (0,T)$.


\subsection{Estimates independent of $R$}\label{limr}

We test \eqref{PDE3} by $w = u_{t}$, \eqref{PDE4} by $z=1$, and sum up.
Unlike in \eqref{cons}, we obtain the inequality 
\begin{eqnarray} \nonumber
&&\frac{\dd}{\dd t} \int_0^\ell \Big(\hat C_V(\theta) + \frac{\rho}{2} u_t^2 + \frac{c}{2} u_x^2
+ \frac{\gamma}{2} u_{xx}^2 + f(u_x)\UTh
+ \frac12 \alpha(u_x)(\PTh)^2\\ \label{ener}
&&\quad + \frac{1}{2\kappa}\big(e u_x +\PTh\big)^2 + \htt f(u_x)(q \PTTh - \UTTh)\Big) \dd x \le 0,
\end{eqnarray}
where
\be{hatc}
\hat C_V(\theta) := \int_0^\theta C_V(K_R(\vartheta))\dd\vartheta \ge c_0\Big(\theta + \frac34\htt^{4/3}\Big)
\ee
by virtue of Hypothesis \ref{h1}\,(iii). Using \eqref{he1}, \eqref{he2}, and \eqref{he7} we have $|\PTh| \le C$,
$0 \le \UTh \le C(1+\htt)$, and $\htt f(u_x)|(q \PTTh - \UTTh)| \le C$, hence
\begin{eqnarray*}
&&\hspace{-10mm}\il \left(\theta + \htt^{4/3} + u_t^2 + u_x^2 + u_{xx}^2\right)(x,t)\dd x
\\
&&\le  C\left(1+\int_0^\ell \Big(\hat C_V(\theta) + \frac{\rho}{2} u_t^2 + \big(\frac{c}{2}+\frac{e}{2\kappa}\big) u_x^2
+ \frac{\gamma}{2} u_{xx}^2\Big)(x,0)\dd x\right),
\end{eqnarray*}
and consequently
\be{ener2}
\il \left(\theta + \htt^{4/3} + u_t^2 + u_x^2 + u_{xx}^2\right)(x,t)\dd x \le C
\ee
for every $t \in [0,T]$.

We now introduce a more convenient notation. For $p\ge 1$ and $v \in L^p(0,\ell)$ we denote
$$
|v|_p = \left(\il |v(x)|^p\dd x\right)^{1/p}.
$$
Similarly, for $v \in L^p((0,\ell)\times (0,T))$ and $\tau \in [0,T]$ we put 
$$
\|v\|_{p,\tau} = \left(\int_0^\tau\!\il |v(x,t)|^p\dd x\dd t\right)^{1/p}.
$$
We can rewrite the inequalities \eqref{es3b}--\eqref{es3c} and \eqref{ener2} in the form
\begin{eqnarray}\nonumber
\hspace{-15mm}|\theta_x(\tau)|_2^2 + |u_{tt}(\tau)|_2^2 + |u_{xt}(\tau)|_2^2 + |u_{xxt}(\tau)|_2^2 && \\ \label{ee1}
+ \int_0^\tau\!\!\il (1+\htt^{1/3})|\theta_{t}|^2 \dd x \dd t + \|u_{xtt}\|_{2,\tau}^2
&\le& C\left(1 {+} \|\htt\|_{4,\tau}^4 {+} \|u_{xt}\|_{4,\tau}^4\right),\\ \label{ee2}
|u_t(\tau)|_2^2 {+} |u_x(\tau)|_2^2 {+} |u_{xx}(\tau)|_2^2
+ \|u_{xt}\|_{2,\tau}^2 &\le& C\left(1 {+} \|\htt\|_{2,\tau}^2\right),\\ \label{ee3}
|\theta(\tau)|_1 + |\htt(\tau)|_{4/3}^{4/3} + |u_t(\tau)|_2^2 + |u_x(\tau)|_2^2 + |u_{xx}(\tau)|_2^2
&\le& C
\end{eqnarray}
for every $\tau \in [0,T]$ with a constant $C>0$ independent of $R$ and $\tau$.

We now repeatedly use the Gagliardo-Nirenberg inequality (see, e.\,g., \cite{bin} for a general information)
in its simplest form, which states that there exists a constant $C>0$ such that for every
$v \in W^{1,p}(0,\ell)$ with $p\ge 1$ and every $1 \le s < q$ we have
\be{gn}
|v|_q \le C(|v|_s + |v|_s^{1-a} |v_x|_p^a)\,, \ a = \frac{\frac1s - \frac1q}{1+ \frac1s - \frac1p} \in (0,1).
\ee
We first observe that for all $\tau$ we have, by virtue of \eqref{ee3}, that
\be{gn1}
|\htt(t)|_2 \le C\left(|\htt(t)|_{4/3} + |\htt(t)|_{4/3}^{4/5}|\htt_x(t)|_{2}^{1/5}\right)
\le C\left(1+|\theta_x(t)|_{2}^{1/5}\right),
\ee
hence,
\be{gn2}
\|\htt\|_{2,\tau}^2 \le C(1+ \|\theta_x\|_{2,\tau}^{2/5}).
\ee
Similarly,
\be{gn3}
|\htt(t)|_4 \le C\left(|\htt(t)|_{4/3} + |\htt(t)|_{4/3}^{3/5}|\htt_x(t)|_{2}^{2/5}\right)
\le C\left(1+|\theta_x(t)|_{2}^{2/5}\right),
\ee
hence,
\be{gn4}
\|\htt\|_{4,\tau}^4 \le C(1+ \|\theta_x\|_{2,\tau}^{8/5}).
\ee
The $\theta_x$-term on the left-hand side of \eqref{ee1} is therefore dominant, and the remaining
critical inequalities read
\begin{eqnarray} \label{ee4}
|\theta_x(\tau)|_2^2 + \|\theta_{t}\|_{2,\tau}^2 + |u_{xxt}(\tau)|_2^2 + \|u_{xtt}\|_{2,\tau}^2
&\le& C\left(1 {+} \|\theta_x\|_{2,\tau}^{8/5}+\|u_{xt}\|_{4,\tau}^4\right),\\ \label{ee5}
\|u_{xt}\|_{2,\tau}^2 &\le& C\left(1 {+} \|\theta_x\|_{2,\tau}^{2/5}\right)
\end{eqnarray}
for every $\tau\in [0,T]$. For $t \in [0,T]$ we have by \eqref{gn}
\be{ee6}
|u_{xt}(t)|_4 \le C\left(|u_{xt}(t)|_2 + |u_{xt}(t)|_2^{3/4}|u_{xxt}(t)|_2^{1/4}\right).
\ee
To estimate the right-hand side of \eqref{ee6}, put
\be{v}
v(t) = |u_{xt}(t)|_2^2.
\ee
For functions $v: [0,T] \to \real$ and numbers $p\ge 1$ we introduce the seminorms
\be{lp}
|v|_{p,\tau} = \left(\int_0^\tau |v(t)|^p\dd t\right)^{1/p}, \quad \tau \in [0,T].
\ee
For each $v \in L^p(0,T)$, $\tau \in (0,T]$, and $t \in (0,\tau)$ we have
$$
|v(t)|^p \le |v(t)| \left(|v(0)| + \supess_{s\in(0,\tau)}|v(s) - v(0)|\right)^{p-1}
\le |v(t)| \left(|v(0)| + \int_0^\tau |\dot v(t)|\dd t\right)^{p-1},
$$
hence
\be{gn5}
|v|_{p,\tau} \le |v|_{1,\tau}^{1/p}(|v(0)| + |\dot v|_{1,\tau})^{1/p'}, \quad \frac1p + \frac1{p'} = 1.
\ee
This is indeed a variant of the Gagliardo-Nirenberg inequality \eqref{gn} for the case that the initial condition is under control.

With the choice \eqref{v}, we can rewrite \eqref{ee6} in the form
\be{ee7}
|u_{xt}(t)|_4^4 \le C\left(|v(t)|^2 + |v(t)|^{3/2}|u_{xxt}(t)|_2\right),
\ee
hence,
\be{ee8}
\|u_{xt}\|_{4,\tau}^4 \le C\left(|v|_{2,\tau}^2 + |v|_{3/2,\tau}^{3/2} \supess_{t\in(0,\tau)}|u_{xxt}(t)|_2\right).
\ee
Using \eqref{gn5} in \eqref{ee8} successively for $p=2$ and $p=3/2$ and considering the fact that
$$
|\dot v|_{1,\tau} = 2\int_0^\tau\il |u_{xt}(x,t) u_{xtt}(x,t)|\dd x\dd t \le 2 \|u_{xt}\|_{2,\tau}\|u_{xtt}\|_{2,\tau},
$$
we obtain
\be{ee9}
|v|_{2,\tau}^2 \le C|v|_{1,\tau}(1+ |\dot v|_{1,\tau}) \le C\|u_{xt}\|_{2,\tau}^2(1+ \|u_{xt}\|_{2,\tau}\|u_{xtt}\|_{2,\tau}).
\ee
{}From \eqref{ee5} it follows that
\be{ee10}
|v|_{2,\tau}^2 \le C(1+\|\theta_x\|_{2,\tau}^{3/5}\|u_{xtt}\|_{2,\tau}).
\ee
Similarly as in \eqref{ee9} we have
\be{ee11}
|v|_{3/2,\tau}^{3/2} \le C|v|_{1,\tau}(1+ |\dot v|_{1,\tau})^{1/2}
\le C\|u_{xt}\|_{2,\tau}^2(1+ \|u_{xt}\|_{2,\tau}\|u_{xtt}\|_{2,\tau})^{1/2},
\ee
so that
\be{ee12}
|v|_{3/2,\tau}^{3/2} \le C(1+\|\theta_x\|_{2,\tau}^{1/2}\|u_{xtt}\|_{2,\tau}^{1/2}).
\ee
It follows from \eqref{ee4}, \eqref{ee8}, \eqref{ee10}, \eqref{ee12} that for every $\tau\in [0,T]$ we have
\begin{eqnarray} \nonumber
&&\hspace{-12mm}\supess_{t\in(0,\tau)}|\theta_x(t)|_2^2 + \|\theta_{t}\|_{2,\tau}^2
+ \supess_{t\in(0,\tau)}|u_{xxt}(t)|_2^2 + \|u_{xtt}\|_{2,\tau}^2
\le C\left(1 {+} \|\theta_x\|_{2,\tau}^{8/5}+\|u_{xt}\|_{4,\tau}^4\right)\\ \label{ee13}
&\le& C\big(1+ \|\theta_x\|_{2,\tau}^{8/5}+
\|\theta_x\|_{2,\tau}^{3/5}\|u_{xtt}\|_{2,\tau}
+ \|\theta_x\|_{2,\tau}^{1/2}\|u_{xtt}\|_{2,\tau}^{1/2}\supess_{t\in(0,\tau)}|u_{xxt}(t)|_2\big).
\end{eqnarray}
The elementary Young inequality $a^{1/2}b^{1/2}c \le C a^2 + \eta(b^2+c^2)$ with a suitably small $\eta$ enables us to reduce
the inequality \eqref{ee13} to
\be{ee14}
\supess_{t\in(0,\tau)}|\theta_x(t)|_2^2 + \|\theta_{t}\|_{2,\tau}^2
+ \supess_{t\in(0,\tau)}|u_{xxt}(t)|_2^2 + \|u_{xtt}\|_{2,\tau}^2
\le C\left(1 {+} \|\theta_x\|_{2,\tau}^{2}\right),
\ee
and the Gronwall argument yields
\be{ee15}
\supess_{t\in(0,\tau)}|\theta_x(t)|_2^2 + \|\theta_{t}\|_{2,\tau}^2
+ \supess_{t\in(0,\tau)}|u_{xxt}(t)|_2^2 + \|u_{xtt}\|_{2,\tau}^2 \le C.
\ee
Standard embedding theorems imply that $\theta, u_{xt} \in C([0,\ell]\times[0,T])$, and their sup-norm is bounded
by a constant $C_*$ independent of $R$. If we choose $R > 1+ C_*^2$, 
we see that the cut-off functions in \eqref{PDE3}--\eqref{PDE4}
are not active and the solution of \eqref{PDE3}--\eqref{PDE4} that we have constructed by Galerkin approximations
is a solution of \eqref{PDE1}--\eqref{PDE2}, too. This concludes the proof of Theorem~\ref{t1}.


\section{Inversion of time-dependent Preisach operators}\label{inv}

Let $\Psi: \R^L\times\R_+\times\R\to\R$ be a measurable function.
The parameter-dependent Preisach operator $\mathcal{P}:C([0,T];\real^L)\times C[0,T]\to C[0,T]$ with density
$\Psi(\tht,r,v)$  is defined by the integral formula 
\begin{equation}\label{P}
\mathcal{P}(\tht)[q](t)=\int_0^\infty\int_0^{\xi_r(t)} \Psi(\tht(t),r,v)\,\dd v\,\dd r,
\end{equation}
for $\tht\in C([0,T];\real^L)$ and $q\in C[0,T]$, where $\xi_r=\play_r[q]$ stands for the solution
of the variational inequality \eqref{vari}. 

To ensure that the definition \eqref{P} is meaningful, it is convenient to reduce the set of admissible functions $\Psi$ by assuming the following:

\begin{hypothesis}\label{basic} For all $\vtht \in \real^L$ we have
\begin{enumerate}
\item[$(i)$] $0\leq \Psi(\vtht,r,v)\le \mu(r)$ a.e., where $\mu\in L^1(0,\infty)$ and $M=\int_0^\infty\mu(r)\dd r$;
\item[$(ii)$] $\left\|\nabla_\vtht\Psi(\vtht,r,v)\right\|\le K(r,v)$ a.e., where 
$K\in L^1((0,\infty)\times\R)$, $\|\cdot\|$ denotes the Euclidean norm in $\real^L$, and 
\[
M_1=\int_0^\infty\int_{-\infty}^\infty K(r,v)\dd v\,\dd r.
\]
\end{enumerate}
\end{hypothesis}
Under these conditions, one can show that the operator $\mathcal{P}$ is Lipschitz continuous (see \cite{km1}). More specifically, we endow the space $C[0,T]$ with a family of seminorms 
$$|w|_{[s,t]} := \sup_{\tau \in [s,t]} |w(\tau)|.$$ 
Then, for $q_1,q_2\in C[0,T]$, $\tht_1,\tht_2\in C([0,T];\real^L)$, and $t\in[0,T]$, we have
\begin{equation*}\label{P-Lip}
\big|\mathcal{P}(\tht_1)[q_1](t)-\mathcal{P}(\tht_2)[q_2](t)\big|
	\leq M\,\mot{q_1-q_2}+ M_1\,\|\tht_1(t)-\tht_2(t)\|
\end{equation*}

It is worth mentioning that in \cite{km1} the theory of parameter-dependent Preisach operator is developed in a more general setting of regulated functions (i.e., functions having only discontinuities of the first kind).

Next, we recall the inversion formula proved in \cite{km1}. 

\begin{theorem}\label{p37}
Assume that a measurable function $\Psi:\R^L\times\R_+\times\R\to\R$ satisfies Hypothesis \ref{basic}.
Let $\tht,\,\dehat{\tht}\in C([0,T];\real^L)$, $w,\,\dehat{w}\in C[0,T]$ be given.
Then there exist solutions $q,\,\dehat{q}\in C[0,T]$ of the equations
\begin{equation}\label{eq77}
\begin{split}
q(t)+\mathcal{P}(\tht)[q](t)=w(t),
\\[2mm]
\dehat{q}(t)+\mathcal{P}(\dehat{\tht})[\dehat{q}](t)=\dehat{w}(t),
\end{split}
\qquad t\in[0,T],
\end{equation}
and the inequality
\begin{equation}\label{eq78}
|q(t)-\dehat{q}(t)|\le \expe^M\,\big(\mot{w-\dehat{w}} + M_1\,\Mot{\tht-\dehat{\tht}}\big)
\end{equation}
holds for each $t\in[0,T]$. 
If, in addition, $\tht$ and $w$ are absolutely continuous,
then the solution $q$ is absolutely continuous.
\end{theorem}

Theorem \ref{p37} deals precisely with the situation in equation \eqref{q}, where $\tht = (\eps,\theta)$, 
\[
w = \frac{e\eps}{\kappa f(\eps)}, 
\]
and
\[
\Psi(\tht,r,v) = \frac{1+ \kappa\alpha(\eps)}{\kappa f(\eps)} 
\psi(\theta, r,v)\,\dd v\,\dd r.
\]
In this case, Hypotheses \ref{h1} and \ref{h2} guarantee that Hypotheses \ref{basic} are satisfied.

\end{document}